\documentclass[12pt]{amsart}
\usepackage{amssymb}
\usepackage{times}
\usepackage{graphicx}

\theoremstyle{plain}

\numberwithin{equation}{section}

\newcommand{\calC}{\mathcal{C}}

\newcommand{\bbF}{\mathbb{F}}

\newcommand{\bbQ}{\mathbb{Q}}
\newcommand{\bbR}{\mathbb{R}}

\newcommand{\bbZ}{\mathbb{Z}}

\def\Aut{{\text{Aut}}}

\def\PSL{{\text{PSL}}}
\def\rank{{\text{rank}}}
\def\det{{\text{det}}}
\def\mod{{\text{mod}}}

\begin{document}
\title [Supersingular K3 surfaces] {Maximal subgroups of the Mathieu group $M_{23}$ and 
symplectic automorphisms of supersingular K3 surfaces}
\author{Shigeyuki Kond{$\bar{\rm o}$}}
\address{Graduate School of Mathematics, Nagoya University, Nagoya,
464-8602, Japan}
\email{kondo@math.nagoya-u.ac.jp}
\thanks{Research of the author is partially supported by
Grant-in-Aid for Scientific Research A-14204001and Hoga-17654004, Japan}

\begin{abstract}
We show that the Mathieu groups $M_{22}$ and $M_{11}$  
can act on the supersingular $K3$ surface
with Artin invariant 1 in characteristic 11 as symplectic automorphisms.   More generally we show that
all maximal subgroups of the Mathieu group $M_{23}$ with three orbits on 24 letters act on a
supersingular $K3$ surface with Artin invariant 1 in a suitable characteristic.
\end{abstract}
\maketitle

\section{Introduction}
Let $X$ be a $K3$ surface defined over an algebraically closed field. By definition, the irregularity of $X$ 
vanishes and there exists a unique (up to constants) non-zero regular 2-form on $X$.
An automorphism $g$ of $X$ is called {\it symplectic} if $g$ fixes a non-zero
regular 2-form on $X$.  In case of complex $K3$ surfaces, Mukai \cite{M} showed that
any finite group of sympletic automorphisms of a $K3$ surface is a subgroup of the Mathieu group $M_{23}$
with at least five orbits in its natural action on 24 letters.  However in case of positive characteristic,
this does not hold.  For example, the projective unitary group $PU(4, \bbF_9)$ acts on 
the Fermat quartic surface in characteristic 3 as projective transformations. 
By comparing their orders we can see that the group $PU(4,\bbF_9)$ is not a subgroup of $M_{23}$.  Note that
the Fermat quartic surface in characteristic 3 is a supersingular $K3$ surface with Artin invariant 1 (Shioda \cite{S}).
Also Dolgachev and the author \cite{DKo} proved that the group $L_3(4) : 2$ acts on a supersingular $K3$ surface
with Artin invariant 1 in characteristic 2.  In this case $L_3(4):2$ is not a subgroup of $M_{23}$, too.
Recently Dolgachev and Keum \cite{DKe1}, \cite{DKe2} studied the details in case of positive characteristic.
In particular they are trying to extend Mukai's theorem to the case of positive characteristic.

In this note, inspired by Dolgachev and Keum \cite{DKe2}, 
we shall show that each maximal subgroup of $M_{23}$ with three orbits on 24 letters 
can act as automorphisms on a supersingular $K3$ surface with
Artin invariant 1 by using Ogus's Torelli type theorem for supersingular $K3$ surfaces (Ogus \cite{O1}, \cite{O2})  (see Theorem \ref{Main}).
The simpleness of $M_{22}$ and $M_{11}$ imply that these actions are symplectic (Corollary \ref{Cor}).
The idea of the proof comes from Mukai's one in the appendix of \cite{K}.
Let $N$ be the Niemeier lattice with root sublattice $A_1^{24}$.  
Here we consider the negative definite one as $N$.  The Mathieu group
$M_{24}$ naturally acts on the set of 24 positive roots of $A_1^{24}$ as permutations and $M_{23}$ is the stabilizer of a fixed positive root.
Let $G$ be a maximal subgroup of $M_{23}$ with 3 orbits on 24 letters.  We can consider $G$ as a subgroup
of the orthogonal group of $N$.  Let $N^G$ be
the invariant sublattice.  Then by assumption $N^G$ is of rank 3, and 
hence the orthogonal complement $N_G$ of
$N^G$ in $N$ is an even negative definite lattice of rank 21 and $N_G$ contains no $(-2)$-vectors.
We can see that there exists an even positive definite lattice $<h>$ of rank 1 such that $<h> \oplus N_G$ can be embedded into the N{\' e}ron-Severi lattice $S_X$ 
of a supersingular $K3$ surface $X$ with Artin invariant 1 in a suitable characteristic $p$.
The action of $G$ on $N_G$ can be extended to the one on $S_X$ acting trivially on $<h>$.  Since $N_G$ contains no $(-2)$-vectors, we may assume that $G$ preserves the ample cone of $X$.  Moreover
$G$ acts trivially on $N_G^*/N_G \cong (N^G)^*/N^G$ and hence acts trivially on $S_X^*/S_X$.
This implies that $G$ preserves the "period" of $X$.  Therefore it follows from Ogus's Torelli
theorem \cite{O2} that $G$ is realized as a subgroup of $\Aut(X)$.

We use the following symbols of finite groups in this paper:
\medskip

$n$ : a cyclic group of order {\it n}.
\medskip

$n^k$ : an $n$-elementary abelian group of order $n^k$.
\medskip

$S_n$  ($A_n$) : a symmetric (alternating) group of degree $n$.
\medskip

$L_{n}(q)$ : the projective special linear group $\PSL(n,q)$.
\medskip

$M_k$ $(k = 11, 12, 22, 23, 24)$  : the Mathieu group.

\medskip
\noindent
We shall say that a group $G$ is a group $N{\cdot}H$ when we mean that $G$ has a normal
subgroup $N$ whose quotient is isomorphic to $H$.  We denote by $N:H$ the semi-direct product.

\medskip
\noindent
{\bf Acknowledgments.}  The author thanks to JongHae Keum for stimulating discussions.

\section{Lattices}\label{}

\subsection{Preliminaries}\label{}
A {\it lattice} is a a free $\bbZ$-module $L$ of finite rank endowed with a $\bbZ$-valued  
symmetric bilinear form $\langle , \rangle$.  If $L_1$ and $L_2$ are lattices, then $L_1\oplus L_2$
denotes the orthogonal direct sum of $L_1$ and $L_2$.  Also we denote by $L^m$ the 
orthogonal direct sum of {\it m}-copies of $L$.   An isomorphism of lattices preserving the
bilinear forms is called an {\it isometry}. For a lattice $L$, we denote by $O(L)$ the 
group of self-isometries of $L$.  A sublattice $S$ of $L$ is called {\it primitive} if
$L/S$ is torsion free.

A lattice $L$ is {\it even} if $\langle x,x \rangle$ is even for each $x\in L$.  A lattice $L$ is 
{\it non-degenerate} if the discriminant $d(L)$ of its bilinear form is non zero, and 
{\it unimodular} if $d(L) = \pm 1$. If $L$ is a non-degenerate lattice, the {\it signature} of $L$ is a pair 
$(t_+, t_-)$ where $t_{\pm}$ denotes the multiplicity of the eigenvalues ${\pm 1}$ 
for the quadratic form on $L \otimes \bbR$.

Let $L$ be a non-degenerate even lattice.  The bilinear form of $L$ determines a canonical embedding $L \subset L^{\ast}={\rm Hom}(L,\bbZ )$. The factor group $L^{\ast}/L$, which is denoted by $A_L$, 
is an abelian group of order $|d(L)|$. We denote by $l(L)$ the
number of minimal generators of $A_L$. We extend the bilinear form on $L$ to the one on
$L^{\ast}$, taking value in $\bbQ$, and define 
$$q_L : A_L \to \bbQ/2 \bbZ,\,\   q_L(x+L)\, =\, \langle x,x \rangle +\,2 \bbZ\,\   (x\in L^{\ast}).$$  
We call $q_L$ the {\it discriminant quadratic form} of $L$.

Let $S$ be an even lattice.  Let $L$ be an even lattice containing $S$ as a sublattice of finite index.
We call $L$ an {\it overlattice }  of $S$.  Note that $L$ is determined by the isotropic subgroup
$L/S$ in $A_S$ with respect to $q_S$.

We denote by $U$ the even lattice defined by the matrix
$
\begin{pmatrix}0&1
\\1&0
\end{pmatrix}
$
and by $A_m$, $D_n$ or $E_l$ the even negative definite lattice defined by the Dynkin matrix of
type $A_m$, $D_n$ or $E_l$ respectively.

\subsection{The N{\' e}ron-Severi lattice of a supersingular $K3$ surface}\label{}

A {\it supersingular} $K3$ surface is a $K3$ surface with the Picard number 22.
A supersingular $K3$ surface exists only in positive charavteristic $p$.
Let $X$ be a supersingular $K3$ surface in characteristic $p$ and let $S_X$ be the
N{\' e}ron-Severi lattice of $X$.  It is known that $\det (S_X) = -p^{2\sigma}$, 
$(1 \leq \sigma \leq 10)$ where the number $\sigma$ is called {\it Artin invariant} of $X$ (Artin \cite{A}).
A generic supersingular $K3$ surface has  Artin invariant 10 and a supersingular $K3$ surface with
$\sigma = 1$  is unique.  Moreover the N{\' e}ron-Severi lattice $S = S_X$ is uniquely determined by $\sigma$
(Rudakov-Shafarevich \cite{RS}, Ogus \cite{O1}).
For example, 
$$S = U \oplus E_8 \oplus A_6 \oplus A_6 \ (p = 7, \sigma = 1);$$
$$S = U \oplus A_{10} \oplus A_{10} \ (p = 11, \sigma = 1).$$

In case $p = 5$ and $\sigma = 1$, $S$ is obtained as follows.  Let $K = U \oplus E_7 \oplus  A_4 \oplus A_9$.
Then $K^*/K \cong (\bbZ / 5\bbZ)^2 \oplus (\bbZ / 2\bbZ)^2$.  Let $x$ be a generator of $E_7^*/E_7$
and $y$ a generator of $A_{9}^*/A_9$.  Then $q_S(x) = 1/2$ and  $q_S(5y) = -1/2$.  
The isotropic vector $x + y$ of $K^*/K$ determines an even lattice $S$ which contains $K$ of index 2 and is
with $\det (S) = -5^2$.

The discriminant form of the above $S$ is as follows:
$$(A_S, q_S) = ((\bbZ /5\bbZ)^2,  (-2/5) \oplus (-6/5)) \ (p = 5, \sigma = 1);$$
$$(A_S, q_S) = ((\bbZ /7\bbZ)^2, (-6/7) \oplus (-6/7) \ (p = 7, \sigma = 1);$$
$$(A_S, q_S) = ((\bbZ /11 \bbZ)^2, (-10/11)  \oplus (-10/11)  \ (p = 11, \sigma = 1).$$

\subsection{Niemeier lattices and Mathieu groups}\label{}

A {\it Niemeier lattice} is an even negative definite unimodular lattice of rank 24.
The isomorphism class of a Niemeier lattice is determined by the sublattice $R$ generated by
all $(-2)$-vectors in it.  It is known that there exists a Niemeier lattice $N$ with $R = A_1^{24}$.
Moreover the orthogonal group $O(N)$ is isomorphic to $2^{24} : M_{24}$.  The subgroup $2^{24}$
is generated by reflections associated to 24 positive roots in  $A_1^{24}$ and 
$M_{24}$ naturally acts on the set of 24 positive roots
of $A_1^{24}$.  Then $M_{23}$ is the stabilizer of a fixed positive root.  The following is the table of
all maximal subgroups of $M_{23}$ (\cite{C}, page 71,  \cite{CS}, Chap. 10).

\begin{table}[h]
\[
\begin{array}{rllll}
{}& {\rm Maximal \ subgroup}&{\rm Order}&{\rm Orbit \ Decomposition} \\
1) & M_{22}&2^7\cdot 3^2 \cdot 5 \cdot 7 \cdot 11 & [1,1,22] \\
\noalign{\smallskip}
2) & L_3(4) : 2 & 2^7\cdot 3^2 \cdot 5 \cdot 7 & [1,2,21] \\
    \noalign{\smallskip}
3) & 2^4 : A_7 & 2^7\cdot 3^2 \cdot 5 \cdot 7 & [1,7,16] \\
    \noalign{\smallskip}
4) & A_8 & 2^6\cdot 3^2 \cdot 5 \cdot 7& [1, 8, 15] \\ 
\noalign{\smallskip}
5) & M_{11}  & 2^4\cdot 3^2 \cdot 5 \cdot 11 & [1,11,12]  \\ 
\noalign{\smallskip}
6) & 2^4 : (3\times A_5) : 2  & 2^7 \cdot 3^2 \cdot 5 & [1, 3, 20] \\
\noalign{\smallskip}
7) & 23 : 11 & 11 \cdot 23 & [1, 23] & \\ 
\end{array}
\]
\caption{}
\end{table}

\subsection{Remark}\label{}
The group $L_3(4):2$ in the Table 1 is different from the one mentioned in Introduction which appeared in
the paper \cite{DKo}.
In the case of Table 1, the involution $2$ is $2_2$ given in \cite{C}, page 71, and in the case of \cite{DKo},
the involution $2$ is $2_1$ given in \cite{C}, page 80.

\medskip

We recall that the Niemeier lattice $N$ is obtained from $A_1^{24}$ as follows.
Let $\calC$ be the binary Golay code which is a subspace of $(A_1^*/A_1)^{24} \cong \bbF_2^{24}$
of dimension 12.   Then 
$$N = \{ x \in (A_1^*)^{24} \ | \ x \ \mod \ A_1^{24} \in \calC \}.$$
It is known that the length of non-zero entries of $x \in \calC$ is 8, 12, 16 or 24.
The set of non-zero entries of length 8 is called an {\it octad} and one of length 12 a {\it dodecad}.
In the case 3 on the Table 1, the union of orbits of length 1 and 7 is an octad.
Also in case 4, the orbit of length 8 is an octad.   In case 5, the orbit of length 12 and its complement are
dodecad.  For more details, we refer the reader to Conway-Sloane \cite{CS}.

\section{Wild symplectic automorphisms}\label{}

In this section we shall prove the following:

\subsection{Theorem}\label{Main}
{\it Let $G$ be a maximal subgroup of $M_{23}$ with three orbits.  Then there exists a prime number $p$ 
such that $G$ acts as automorphisms on 
a supersingular $K3$ surface with Artin invariant $1$ in characteristic $p$.}

\medskip

First we shall show the following Lemma.  

\subsection{Lemma}\label{Key}
{\it Let $G$ be a maximal subgroup of $M_{23}$ with three orbits.  Then there exists a prime number $p$
such that $G$ acts on the N{\' e}ron-Severi lattice $S$ of a supersingular $K3$ surface 
with Artin invariant $1$ in characteristic $p$.  Moreover $G$ acts trivially on $S^*/S$ and the orthogonal
complement of the invariant sublattice $S^G$ in $S$ contains no $(-2)$-vectors.}

\begin{proof}
Let $N$ be the Niemeier lattice with the root sublattice $A_1^{24}$ on which $M_{23}$
naturally acts.  Let $N^G$ be the invariant sublattice.  Since $G$ has three orbits,
$\rank (N^G) = 3$.  Let $N_G$ be the orthogonal complement of $N^G$ in $N$.  Then
$\rank (N_G) = 21$.  For each $G$ in the Table 1, we shall show the following:
First we calculate the discriminant forms $q_{N_G} = -q_{N^G}$.
Next we take a vector $h$ with $h^2 = |\det (N_G) |$ and consider the lattice 
$<h> \oplus N_G$.  Then we shall show that there exists an over lattice $S$ of $<h> \oplus N_G$ which is isomorphic to
the N{\' e}ron-Severi lattice of a supersingular $K3$ surface $X$.  Moreover
the action of $G$ on $N_G$ can be extended to the one on $S$ acting trivially on $<h>$.
Since $G$ acts on $<h>^*/<h>  \oplus N_G^*/N_G$ trivially, $G$ acts on $S^*/S$ trivially.  

Note that $N$ contains exactly 24 positive roots ($(-2)$-vectors) and $G$ acts on the set of positive roots as permutations.
Hence $N_G$ contains no $(-2)$-vectors.

In the following we denote by $\{ x_1,..., x_{24} \}$ the set of positive roots of $A_1^{24}$.
\smallskip

{\it Case} 1:  $G = M_{22}$.   

We assume that $x_1, x_2, x_3 +  \cdot \cdot \cdot + x_{24}$ are invariant
under the action of $G$.   Then $N^G$ is generated by $x_1, x_2, x_3 +  \cdot \cdot \cdot + x_{24}$ and
$(x_1 + x_2 + x_3 + \cdot \cdot \cdot + x_{24})/2$.  Hence $\det (N^G) = 2\cdot 2 \cdot 44/2^2 = 44$.
By using Nikulin \cite{N1}, Proposition 1.5.1, we can easily see that  $q_{N^G} = (-5/4) \oplus (-4/11)$.  
Hence $q_{N_G} =  -q_{N^G} = (5/4) \oplus (4/11)$ (Nikulin \cite{N1}, Corollary 1.6.2).
Take a vector $h$ with $h^2 = 44$.  Consider the subgroup $H$ of order 4 in
$<h>^*/<h>  \oplus N_G^*/N_G$ generated by $h/4 + \theta$, where $\theta$ is a generator of
2-Sylow subgroup of $N_G^*/N_G$.
Since $H$ is totally isotropic with respect to $q_{<h>}  \oplus q_{N_G}$, it determines the overlattice $S$
with the discriminant form $q_S =  (-10/11)  \oplus (-10/11)$ (Nikulin \cite{N1}, Proposition 1.4.1).  
It now follows from Nikulin \cite{N1}, Theorem 1.14.2 that
$S$ is isomorphic to the N{\' e}ron-Severi lattice of the supersingular $K3$ surface with Artin invariant 1 in
characteristic 11.
\smallskip

{\it Case} 2:  $G = L_3(4) : 2$. 

We assume that $x_1, x_2+x_3, x_4+\cdot \cdot \cdot + x_{24}$ are invariant
under the action of $G$.   Then $N^G$ is generated by $x_1, x_2+x_3, x_4 +  \cdot \cdot \cdot + x_{24}$ and
$(x_1 + x_2 + x_3 + \cdot \cdot \cdot + x_{24})/2$.  Hence $\det (N^G) = 2\cdot 2^2 \cdot 42/2^2 = 84$
and $q_{N^G} = (-3/4) \oplus (-2/3) \oplus (-6/7)$.  Hence
$q_{N_G} = (3/4) \oplus (2/3) \oplus (6/7)$.  Take a vector $h$ with $h^2 = 84$.
We consider the totally isotropic subspace $H$ of order 12 generated by $h/12 + \theta$ where
$\theta$ is a generator of the subgroup of order 12 in $N_G^*/N_G$.
Then as in the Case 1, $H$ determines the overlattice $S$ isomorphic to the N{\' e}ron-Severi lattice of
a supersingular $K3$ surface with Artin invariant 1 in characteristic 7.
\smallskip

{\it Case} 3:  $G = 2^4 : A_7$.

We assume that 
$x_1, x_2 +\cdot \cdot \cdot + x_8 , x_9 +  \cdot \cdot \cdot + x_{24}$ are invariant
under the action of $G$.   Then $N^G$ is generated by $x_1, x_2 + \cdot \cdot \cdot + x_8, x_9 +  \cdot \cdot \cdot + x_{24}$, 
$(x_1 + x_2 + x_3 + \cdot \cdot \cdot + x_{24})/2$ and  $(x_1 + \cdot \cdot \cdot + x_8)/2$.  
Hence $\det (N^G) = 2\cdot 14 \cdot 32/2^4 = 56$
and $q_{N^G} = (-1/8) \oplus (-2/7)$.  Hence
$q_{N_G} =  (1/8) \oplus (2/7)$.
Take a vector $h$ with $h^2 = 56$.
We consider the totally isotropic subspace $H$ of order 8 generated by $h/8 + \theta$ where
$\theta$ is a generator of the 2-Sylow subgroup of order 8 in $N_G^*/N_G$.
Then as in the Case 1, $H$ determines the overlattice $S$ isomorphic to the N{\' e}ron-Severi lattice of
a supersingular $K3$ surface with Artin invariant 1 in characteristic 7.
\smallskip

{\it Case} 4:  $G= A_8$. 

We assume that 
$x_1, x_2 +\cdot \cdot \cdot + x_9 , x_{10} +  \cdot \cdot \cdot + x_{24}$ are invariant
under the action of $G$.   Then $N^G$ is generated by $x_1, x_2 + \cdot \cdot \cdot + x_9, 
x_{10} +  \cdot \cdot \cdot + x_{24}$, 
$(x_1 + x_2 + x_3 + \cdot \cdot \cdot + x_{24})/2$ and $(x_2 + \cdot \cdot \cdot + x_9)/2$.  
Hence $\det (N^G) = 2\cdot 16 \cdot 30/2^4 = 60$ and
$q_{N^G} = (-1/4) \oplus (-4/3) \oplus (-6/5)$.  Hence
$q_{N_G} = (1/4) \oplus (4/3) \oplus (6/5)$.
Take a vector $h$ with $h^2 = 60$.
We consider the totally isotropic subspace $H$ of order 12 generated by $h/12 + \theta$ where
$\theta$ is a generator of the subgroup of order 12 in $N_G^*/N_G$.
Then as in the Case 1, $H$ determines the overlattice $S$ isomorphic to the N{\' e}ron-Severi lattice of
a supersingular $K3$ surface with Artin invariant 1 in characteristic 5.
\smallskip

{\it Case} 5:  $G = M_{11}$.

We assume that 
$x_1, x_2 +\cdot \cdot \cdot + x_{12} , x_{13} +  \cdot \cdot \cdot + x_{24}$ are invariant
under the action of $G$.   Then $N^G$ is generated by $x_1, x_2 + \cdot \cdot \cdot + x_{12}, 
x_{13} +  \cdot \cdot \cdot + x_{24}$, 
$(x_1 + x_2 + x_3 + \cdot \cdot \cdot + x_{24})/2$ and $(x_1 + \cdot \cdot \cdot + x_{12})/2$.  
Hence $\det (N^G) = 2\cdot 22 \cdot 24/2^4 = 66$ and
$q_{N^G} = (-3/2) \oplus (-2/3) \oplus (-2/11)$.  Hence
$q_{N_G} = (3/2) \oplus (2/3) \oplus (2/11)$.
Take a vector $h$ with $h^2 = 66$.
We consider the totally isotropic subspace $H$ of order 6 generated by $h/6 + \theta$ where
$\theta$ is a generator of the subgroup of order 6 in $N_G^*/N_G$.
Then as in the Case 1, $H$ determines the overlattice $S$ isomorphic to the N{\' e}ron-Severi lattice of
a supersingular $K3$ surface with Artin invariant 1 in characteristic 11.
\smallskip

{\it Case} 6:  $G = 2^4 : (3\times A_5) : 2$.

We assume that $x_1, x_2+x_3 + x_4, 
x_5+\cdot \cdot \cdot + x_{24}$ are invariant
under the action of $G$.   Then $N^G$ is generated by 
$x_1, x_2+x_3 + x_4, x_5 +  \cdot \cdot \cdot + x_{24}$ and
$(x_1 + x_2 + x_3 + \cdot \cdot \cdot + x_{24})/2$.  Hence $\det (N^G) = 2\cdot 6 \cdot 40/2^2 = 120$
and $q_{N^G} = (-9/8) \oplus (-2/3) \oplus (-8/5)$.  Hence
$q_{N_G} = (9/8) \oplus (2/3) \oplus (8/5)$.
Take a vector $h$ with $h^2 = 120$.
We consider the totally isotropic subspace $H$ of order 24 generated by $h/24 + \theta$ where
$\theta$ is a generator of the subgroup of order 24 in $N_G^*/N_G$.
Then as in the Case 1, $H$ determines the overlattice $S$ isomorphic to the N{\' e}ron-Severi lattice of
a supersingular $K3$ surface with Artin invariant 1 in characteristic 5.
\end{proof}

\begin{proof}  (Theorem \ref{Main})
Let $G$ and $S$ be as in Lemma \ref{Key}.  Let $X$ be the supersingular $K3$ surface with
Artin invariant 1 in characteristic $p$ satisfying $S_X \cong S$.  
Since $N_G$ contains no $(-2)$-vectors, $h$ is contained in a
fundamental chamber of the reflection subgroup $W(X)$ of $O(S_X)$ generated by $(-2)$-reflections.
Hence there exists a $w \in W(X)$ so that $w(h)$ is an ample class.  
Thus $wGw^{-1}$ preserves the ample cone of $X$.
Since both $G$ and $W(X)$ act trivially on $S^*/S$, so is $wGw^{-1}$, and 
hence $wGw^{-1}$ preserves the characteristic subspace ("Period") of $X$ (see Ogus \cite{O2}, page 366).  
Now the assertion follows from Ogus \cite{O2}, Corollary of Theorem II' (page 371).
\end{proof}

\subsection{Corollary}\label{Cor}
{\it The Mathieu groups $M_{22}$, $M_{11}$ and the alternating group $A_8$ act as symplectic automorphisms 
on a supersingular $K3$ surface with Artin invariant $1$.}

\begin{proof} 
Since automorphisms act on a regular 2-form on a $K3$ surface as a multiplicative group, the symplecticness follows from
the simpleness of $M_{22}, M_{11}, A_8$.
\end{proof}

We summarize the prime number $p$ and the degree $h^2$ of the invariant polarization $h$ under $G$
in the following Table 2:

\begin{table}[h]
\[
\begin{array}{rllll}
{}& G & p & h^2 \\
1) & M_{22}& 11  & 44 \\
\noalign{\smallskip}
2) & L_3(4) : 2 & 7 & 84 \\
    \noalign{\smallskip}
3) & 2^4 : A_7 & 7 & 56 \\
    \noalign{\smallskip}
4) & A_8 & 5 & 60 \\ 
\noalign{\smallskip}
5) & M_{11}  & 11 & 66  \\ 
\noalign{\smallskip}
6) & 2^4 : (3\times A_5) : 2  & 5 & 120 \\
\end{array}
\]
\caption{}
\end{table}

It would be interesting to realize these actions geometrically.

\subsection{Problem}\label{}  Let $g$ be an automorphism of a $K3$ surface $X$.
In case that $X$ is a complex $K3$ surface, if $g$ is symplectic, then $g$ acts trivially
on the transcendental lattice of $X$ and hence trivially on the discriminant group $S_X^*/S_X$ of
the N{\' e}ron-Severi lattice of $X$ (Nikulin \cite{N2}, Theorem 3.1).   Moreover
if $G$ is a finite group of symplectic automorphisms of a complex $K3$ surface $X$, denote by
$L_G$ the orthogonal complement of the invariant sublattice of $H^2(X, \bbZ)$.
Then 
$$l(L_G) \leq 22 - {\rm rank} (L_G)$$
where $l(L_G)$ is the number of minimal generator of $L_G^*/L_G$ (\cite{K}, Proposition 2).  
This means that if $G$ becomes bigger, then $L_G$ becomes bigger, too, 
and hence $l(L_G)$ becomes smaller.

In case of positive characteristic, does 
any symplectic automorphism of a supersingular $K3$ surface act trivially
on the discriminant group of the N{\' e}ron-Severi lattice ?  And if $G$ becomes bigger, then
does the Artin invariant become smaller ?

\today

\end{document}